\documentclass[10pt]{amsart}
\usepackage{amsmath,amssymb,enumerate}
\usepackage{epsf,epsfig,amsfonts,graphicx,color}
\usepackage[applemac]{inputenc}
\usepackage{url,a4wide,epsfig}
\usepackage{ulem}
 \numberwithin{equation}{section}
\newcommand\op[1]{\mathop{\rm #1}\nolimits}
\newcommand\p{\partial}
\newcommand\R{{\mathbb R}}
\newcommand\ff{{\boldsymbol f}}
\newcommand\bh{{\boldsymbol h}}
\newcommand\bx{{\boldsymbol x}}
\newcommand\by{{\boldsymbol y}}

\newtheorem{theorem}{Theorem}
\newtheorem{remark}[theorem]{Remark}
\newcommand{\weg}[1]{}

\newtheorem{proposition}[theorem]{Proposition}
\newtheorem{cor}[theorem]{Corollary}

\title{Almost every path  structure is not variational}
\author{Boris S.\ Kruglikov and Vladimir S.\ Matveev}
 \address{B.\,S. Kruglikov: Department of Mathematics and Statistics,
UiT the Arctic University of Norway, 9037 Troms\o, Norway.
Email: {\tt boris.kruglikov@uit.no}}
 \address{V.\,S. Matveev: Institute of Mathematics, Friedrich-Schiller-Universit\"at, 07737 Jena, Germany.
Email: {\tt vladimir.matveev@uni-jena.de}}

\begin{document}

 \begin{abstract}
Given a smooth family of unparameterized curves such that through every point in every direction
there passes exactly one curve, does there exist a Lagrangian with extremals being precisely this family?
 It is known that in  dimension 2  the answer is positive.  
In dimension 3, it follows from the work of Douglas that the answer is, in general, negative.
We generalise this result to all  higher  dimensions and show that the answer is actually negative for
almost every such a family of  curves, also known as path structure or path  geometry.

On the other hand, we consider path geometries possessing infinitesimal symmetries and   show that 
 path and projective structures with submaximal symmetry dimensions are variational.  
Note that the  projective structure with the  submaximal symmetry algebra, 
the so-called Egorov structure, is not pseudo-Riemannian metrizable; 
we show that it is metrizable in the class of Kropina pseudo-metrics   and explicitly construct the corresponding Kropina Lagrangian.
 \end{abstract}

\maketitle

\vspace{-0.5cm}

\section{Introduction}\label{S1}

\subsection{ Definitions and  motivations}\label{S11}

Consider the following system of second order ODEs on a space $M$ of dimension $n+1$ with coordinates
$\bx=(x^0,\dots,x^n)$:
 \begin{equation}\label{sODE}
x^i_{tt}+ h^i(\bx,\bx_t)=\nu x^i_t,   \qquad 0\leq i\leq n.
 \end{equation}
Here functions $h^i(\bx,\bx_t)$ are   assumed to be  positively homogeneous of the second degree in $\bx_t$, 
i.e., $\bh(\bx, \lambda \bx_t)= \lambda^2 \bh(\bx,  \bx_t) $ for every $\lambda >0$,  
and $\nu$ is an arbitrary functional parameter to be eliminated.   
That is, a solution of the system is a vector-function $ \bx(t)$ such that there exists a function $\nu(t)$ 
for which \eqref{sODE} holds; $\bx_t$ and $\bx_{tt}$ denote the first and second derivatives
of the vector-function $\bx(t)$ in $t$. This system is clearly underdetermined and effectively consists  
of $n$ equations on $n+1$ unknown functions.  
From the physical viewpoint it can be interpreted as the condition that at every point the acceleration   
$\bx_{tt}+ \bh(\bx,\bx_t)$ is linearly dependent with the velocity $\bx_t$. 

Since $\bh$ is positively homogeneous of the second degree in $\bx_t$, for every solution $\bx(t)$ 
of system \eqref{sODE} and for any local diffeomorphismn $t\mapsto\tau(t)$ of $\R$ with $\tau'(t)>0$ 
the reparameterized curve $\bx(\tau(t))$ is also a solution. 
Therefore, solutions of \eqref{sODE} are arbitrary orientation-preserving reparameterizations of solutions 
of the system  
 \begin{equation}\label{sODEs}
x^i_{tt}+  h^i(\bx,\bx_t)=0,   \qquad 0\leq i\leq n.
 \end{equation}
For any point and any oriented direction there exists exactly one solution with these initial data.

A {\it path structure}  (also known  as {\it path geometry}) is the solution space of an equation of the form \eqref{sODE} or equivalently
of \eqref{sODEs} where we forget parametrization of solution-curves (henceforth called paths).
Geometrically, it is defined as a smooth family of {\it unparameterized} curves such that
there exists precisely one curve from the family through every point in every  oriented  direction.

The simplest example of a path structure is the {\it flat} structure on an affine space,
where all the curves of the family are straight lines. (A locally equivalent path structure is given by 
the geodesic family on a Riemannian space of  constant sectional curvature.)
 
We say that a path structure is {\it reversible}, if for every point and any oriented direction the path passing 
through this point in this direction geometrically coincides with the path passing in the reversed direction. 
For example, the flat structure is reversible. 
Clearly, reversibility is equivalent to the property that for every $(\bx,\bx_t)$ the difference 
$\bh(\bx, \bx_t)- \bh(\bx, -\bx_t)$ is proportional to $\bx_t$.  

Path structures naturally appear in differential geometry and in mathematical relativity.
Indeed, for a Lagrangian\footnote{We use  hat on  autonomous  Lagrangians (which for most part of the paper can be assumed to be homogeneous in $\bx_t$)  
to distinguish them from nonautonomous  Lagrangians  in a space of one dimension less  used later on.} $\hat{L}(\bx,\bx_t)$ 
positively homogeneous of degree one in velocities
(that is $\hat{L}(\bx,\lambda \bx_t)=\lambda\hat{L}(\bx,\bx_t)$ for $\lambda>0$) and  
such that for the ``energy function'' $\hat E:= \tfrac{1}{2} \hat L^2$  the   Hessian $\left(\tfrac{\partial^2  \hat E}{\partial x^i_t \partial x^j_t}\right)$
with respect to $\bx_t$ is nondegenerate, the Euler-Lagrange equation is algebraically-equivalent to a system of the form \eqref{sODE}.
Since unparametrized geodesics of Riemannian and pseudo-Riemannian metrics are extremals of the
Langrangian $\hat{L}$ equal to the square root of the kinetic energy, system \eqref{sODE} contains
the equation of geodesics as a special case. The same is true in Finsler geometry (and pseudo-Finsler
generalisations), where geodesics are extremals of the Lagrangian $\hat{L}$ equal to the Finsler norm;   if the Finsler norm is only positively homogeneous the corresponding path structure can be irreversible. 

Investigation of path structures, as differential equations, and in particular their symmetries, 
goes back to the works of S.~Lie \cite{Lie} and his student A. Tresse \cite{Tresse}. For a scalar ODE of the form 
 \begin{equation} \label{ODE1} 
y_{xx} =  f(x, y, y_x),
 \end{equation}
they considered the path structure on $\mathbb{R}^2(x,y)$ whose paths are given by $x\to (x,y(x))$, 
where $y(x)$ is a solution of \eqref{ODE1}. This path structure is singular in the sense that  in the vertical direction the paths   are   not defined. Symmetries of this  path structure are called {\it point transformations} of the ODE; they correspond to changes of  variables    mixing  dependent and independent variables.

In the context of metric geometry, path structures were studied by H.\ Busemann \cite{Busemann_book}; 
one of the question he considered is whether for a given  path structure there exists 
a Finsler metric whose unparametrized geodesics are paths.

A {\it projective structure} is a path structure given by equation \eqref{sODE} with the functions $h^i$ of the form
 $$
h^i(\bx,\bx_t)= \sum_{j,k=0}^n\Gamma^i_{jk}(\bx)x_t^jx_t^k.
 $$
The corresponding paths are unparametrized geodesics of the  affine connection $(\Gamma^i_{jk})$. Clearly, it is a reversible path structure.
 Projective equivalence of affine
connections is their equivalence as path structures, and was studied since   H.~Weyl \cite{Weyl1} who in particular proved that in dimension $n+1\ge 3$ 
 the Weyl projective curvature tensor $W^\ell_{ijk}$ vanishes if and only if the projective structure is flat.  See also  E.\ Cartan \cite{C}, who constructed the fundamental systems of differential invariants for    projective structures (in dimension $n+1>2$; the case   $n=1$ is due to  \cite{Tresse}).    

A closely related classical problem is when two different metrics have the same geodesics viewed as unparameterized curves.
First nontrivial results in this direction are due to E.\ Beltrami \cite{Beltrami} who proved
that in dimension two a Riemannian metric generating a flat projective structure has constant curvature,
and to U.\ Dini \cite{Dini} who gave a local description of pairs of 2-dimensional Riemannian metrics sharing
the same (unparameterized) geodesics. Results of Beltrami and Dini were generalised to all dimensions by
F.\ Schur \cite{Schur} and T.\ Levi-Civita \cite{LC}.
 % More recently this problem was related to integrability and dynamical behavior of geodesic flows \cite{MT,KM3}.

In the framework of mathematical relativity, projective structures were studied since H.~Weyl \cite{Weyl, Weyl1}. 
He proposed to base the geometric framework of gravity theory on the observable structures of particle trajectories 
and light propagation, i.e., on 
umparameterized geodesics and the conformal structure, see also O.\ Veblen and T.\ Thomas \cite{Veblen}. 
In a fundamental and widely read paper \cite{EPS} J.\ Ehlers, F.\ Pirani and A.\ Schild
claimed that a projective structure and a conformal structure on a differentiable manifold $M$ determine a Weylian metric (Weyl structure),
if and only if the {\it light-like geodesics} of the conformal structure are   paths   of the projective structure.
This claim has been recently proven in \cite{Scholz}; see also \cite{Read},  \cite{matveev} and  \cite[\S 12]{burns}.

Path structures which are not projective structures also naturally appear within mathematical relativity, see the survey by Ch.\ Pfeiffer \cite{Pfeiffer}.
In particular, according to the Fermat principle, projection of null geodesics of a stationary spacetime to a Cauchy hypersurface are geodesics of a Randers (Finsler) metric, see e.g.\ E.\ Caponio et al \cite{Caponio}. These geodesics come without preferred parametrization, since a parameterization depends
on the choice of a Cauchy hypersurface. Note that    path structures coming from most Randers metrics are not reversible; moreover, 
 if a path structure coming from a Randers metric is not reversible, then  
one can uniquely reconstruct this metric up to a trivial projective change   by    \cite{matveev2}.

In our paper we discuss the question whether a given path structure is variational, that is whether there exists a Lagrangian function
{ $\hat{L}(x,x_t, x_{tt},...)$} whose extremals are precisely the paths of the structure. This question is important, because many physical systems can be described
mathematically with the help of the Hamilton-Jacobi formalism and was considered already by H. Helmholz \cite{H}. In differential geometry, this question was explicitly asked by H.\ Busemann \cite{Busemann0}.

In the calculus of variations, this question is one of the so-called {\it inverse problems}, and there is a vast literature on this topic, see e.g.\ books  by I.\ Anderson and G.\ Thompson \cite{Anderson} and  by J.\ Grifone and Z.\ Muzsnay \cite{GM} for two different approaches to this problem (note that the second reference treats 
mostly  parametrized solution-curves of differential equations and is  not directly applicable to our problem), 
as well as the recent surveys \cite{Do1, Do} by T.\ Do and G.\ Prince.

\subsection{Results.}\label{S12}

We consider a path structure in dimension $n+1$ and ask whether there exists an autonomous Lagrangian  such that every curve of our path structure (with any parameterization),  is an extremal of the Lagrangian and  vice versa.  We will call such path structures  {\it variational}.

Our first result shows that we can eliminate higher order derivatives in the Lagrangian $\hat{L}$:

 \begin{proposition}\label{prop:1}
Suppose a path structure  is variational. Then it is variational in the class of Lagrangians of order one:
there exists a positively homogeneous of degree one in velocities function $\hat{L}(\bx,\bx_t)$ whose extremals are precisely curves of the path structure.
 \end{proposition}

Next, we will reduce the problem to a similar one, but in dimension one less. In this reduced problem we will look for nonautonomous Lagrangians
(such a reduction was used in e.g.\ \cite{Beltrami, Lie}, see also \cite[\S 3]{Dunajski}). In order to do this, we parametrize the curves of our path structure
by the first coordinate $x^0=x$ (this is possible locally for almost all solutions). In the notations $\by=(y^1,\dots,y^n)$, $y^j=x^j$ for $1\leq j\leq n$,
the curves are given by $x\mapsto (x,\by(x))$. Thus a path structure on a manifold $M$ is given by a system of second order ODEs,
which in local coordinates can be written as follows (dot means the derivative by $x$):
 \begin{equation}\label{ODE}
\ddot{y}^i=f^i(x,\by,\dot{\by}),\qquad 1\leq i\leq n.
 \end{equation}
   Paths of the path structure are the curves of the form $x\mapsto (x,\by(x))$, where $\by(x)$ are solutions of \eqref{ODE}.  On the language of geometric  theory of ODEs, local diffeomorphisms of the space $(x,\by)$ preserving the path structure  are  called {\it point transformations}.

We will recall in \S\ref{S3} relations between systems \eqref{sODE} and \eqref{ODE}, and explain that
the  inverse variational problem for both systems is  essentially the same. We treat it in the
second (reduced) version. The corresponding Lagrangian $L$ is
a function on the ray-projectivized (or spherical) tangent bundle $STM$.

Recall that the space $J^k(\R^{2n+1},\R^n)$ of $k$-jets of vector-functions
$\ff=(f^i(z))_{i=1}^n$ of the argument $z=(x,\by,\dot{\by})$
consists of the values of independent and dependent variables and their derivatives up to order $k$.
The jet-lift of $\ff$ is the map $j^k\ff:\R^{2n+1}\to J^k$, $z\mapsto(z,\{\partial^j\ff(z)\}_{j=0}^k)$.

 \begin{theorem} \label{Thm}
Let $\ell=4$ for $n>2$ and $\ell=5$ for $n=2$. There exists an open dense set
$\mathcal{U}\subset J^\ell$ such that if $j^\ell\ff(U)\cap\mathcal{U}\neq\emptyset$
for $U\subset STM$ for the right-hand side of \eqref{ODE} then the path   structure  of \eqref{ODE}
is not variational via a first-order Lagrangian $L(x,\by,\dot{\by})$ even microlocally on $U$.
 \end{theorem}

It is well-known that fibers of the bundle $J^k\to J^1$ carry a natural affine structure,
while fibers of $J^1\to J^0$ can be identified with (open charts in) Grassmanians, see e.g.\ \cite{KL}.
Hence, fibers of $J^k\to J^0$ are algebraic, so we can use the Zariski topology.
Recall that open sets in a Zariski topology are open dense in the standard topology, and
the above set $\mathcal{U}$ can be taken Zariski open.
This straightforwardly implies the following statement:

 \begin{cor}\label{Cor} In dimension $n+1$, a generic smooth path structure in $C^4$ topology for $n\ge3$
 and in $C^5$ topology for $n=2$ is not variational (hence not Finsler).
 \end{cor}

In other words, in proper topology, every path structure $\mathcal{P}$ can be  deformed   by an arbitrary small deformation  to a nonvariational
path structure $\mathcal{\tilde P}$ and any sufficiently small deformation  of $\mathcal{\tilde P}$ remains nonvariational.

 % Recall that $C^k$ topology is finer than $C^l$ topology for $k>l$, so
 % for $n\ge3$ a generic smooth path structure in $C^5$ topology is not variational.
 % Thus, for $n>2$ every path structure $\mathcal{P}$ can be deformed to a nonvariational path structure
 % $\mathcal{\tilde P}$ in arbitrary small neighborhood in $C^4$-topology. Moreover, any sufficiently small
 % perturbation of $\mathcal{\tilde P}$ in $C^5$-topology is also nonvariational.

Let us now discuss the dimension $n+1=2$. It is known since 1886, see  
N.~Sonin \cite{Sonin} and G.\ Darboux \cite{Dar}, that in this case every equation \eqref{ODE} 
is (equivalent to) the Euler-Lagrange equation of an nonautonomous Lagrangian. This result was improved in  
\cite{Berck} where it was shown that for  every reversible  path structure $\mathcal{P}$  there exists  a reversible  Finsler metric whose geodesics are paths of the  structure. The irreversible case is still open, see e.g. \cite{Tabachnikov} where the case when all paths are circles was investigated in details.  
 
The case $n+1=3$ was considered by Douglas \cite{D},   who in particular constructed the first example of a nonvariational projective structure.   %with a complete analysis of variationability.
 He also discussed the PDE system
for the inverse variational problem in the case of general $n$, but did not investigate it in detail.
We recall this fundamental system in \S\ref{S3} and in \S\ref{S4} we show how to exploit it to for specific path structures and for all dimensions.

Let us now discuss the question whether all the curves of a given path structure are geodesics of some pseudo-Riemannian metric. In the literature,
this problem is known as ``metrizability'' or . Of course, in this case we may assume that the path structure is actually projective.

Our way to prove Theorem \ref{Thm} easily implies:

 \begin{cor}\label{Cor1}
In dimension $n+1$, a generic smooth projective structure in $C^4$ topology for $n\ge 3$ and
in $C^5$ topology for $n=2$ is not variational, hence not metrizable.
 \end{cor}

The last portion of our results concern path/projective structures with large Lie algebras of symmetries. 
Recall that {\it symmetry} of a path or projective structure is a local diffeomorphism that sends 
paths to paths. It is known that the flat structure in dimension $n+1$ has maximal symmetry dimension 
(i.e., dimension of the symmetry algebra) equal to $n^2+4n+3$. Of course, this path structure is variational 
since geodesics of the Lagrangian $\sqrt{(x_t^0)^2 +...+ (x_t^n)^2\vphantom{\frac{o}{o}}}$ are straight lines.

The next possible symmetry dimensions are $n^2+5$ (for general path structures) and $n^2+4$ 
(for projective  structures), see  \cite{KT}. 
In \S\ref{S51}--\ref{S52} we will demonstrate that these structures are variational 
by exhibiting Lagrangians (of Kropina type). However they are not metrizable: 
for the submaximally symmetric projective structure, called Egorov structure, this follows from  \cite{KM1}; 
the submaximally symmetric path structure is not a projective structure hence can not be metrizable 
by any pseudo-Riemannian metric. This implies the following result: 
  
 \begin{cor}\label{Prop0}
In dimension $n+1>1$ there exists a projective structure that is variational, but not metrizable. 
 \end{cor}

Note that \S \ref{S51} implies this results for $n\ge  2$. For $n=1$, the result is known and follows from e.g.\   
R.\ Bryant et al \cite{Bryant0, Bryant}. Note also that $(n+1=2)$-dimensional   projective  structures admitting  infinitesimal symmetries and the metrization problem  for them   was  solved completely in \cite{Bryant0, alone}. As mentioned above, 2-dimensional projective structures are always variational.   J.\ Lang in  \cite{Lang} constructed Lagrangians for 2-dimensional  path and projective structures with the submaximal symmetry algebra (of dimension $3$), see also \cite{Kr,T}.   

We will also show that the Egorov projective structure is not (regular) Finsler metrizable.
We expect, in the spirit of our results above, that generic variational projective structures are not metrizable
(neither via pseudo-Riemannian nor via Finsler metrics).
We briefly discuss other examples in \S\ref{S53} in relation to the inverse variational problem.

\section{Proof of Proposition \ref{prop:1}}\label{S2}

By the Vainberg-Tonti  formula \cite{Krupka}, if the second order ODE system \eqref{sODE} is variational,
 % (equivalent to the  Euler-Lagrange equation),
then without loss of generality we may assume that the   Lagrangian has the form $\hat{L}=\hat{L}(\bx,\bx_t,\bx_{tt})$.
The corresponding Euler-Lagrange equation then reads:
 \begin{equation}\label{eq:EL}
\frac{d^2}{dt^2} \frac{\p\hat{L}}{\p x^i_{tt}} - \frac{d}{dt} \frac{\p\hat{L}}{\p x^i_t} + \frac{\p\hat{L}}{\p x^i }=0.
 \end{equation}
In this formula the possible highest $t$-derivative of $\bx$ has order $4$ and can come from the terms
$\frac{d^2}{dt^2} \frac{\p\hat{L}}{\p x^i_{tt}}$ only.
Since \eqref{sODE} does not have terms involving  $x^i_{tttt}$, $\hat{L}$ must have the following form:
  \begin{equation}\label{eq:EL1}
\hat{L}(\bx,\bx_t,\bx_{tt})=F(\bx,\bx_t)+\sum_s x_{tt}^s \lambda_s(\bx,\bx_t).
 \end{equation}
Let is now look on the third $t$-derivatives of $\bx$: since the terms with $x^i_{ttt}$ in the equation \eqref{eq:EL}
with $\hat{L}$ given by \eqref{eq:EL} must cancel, we obtain:
 \begin{equation}\label{eq:EL2}
\sum_s \left( \frac{\p\lambda_s }{\p x^i_t} - \frac{\p\lambda_i }{\p x^s_t} \right)x_{tt}^s =0.
  \end{equation}
Then there exists a function $\Lambda(\bx,\bx_t)$ such that $\lambda_s= \tfrac{\p\Lambda}{\p x_t^s}$ implying
$\sum_s x_{tt}^s \tfrac{\p\Lambda}{\p x_t^s}= \tfrac{d}{dt} \Lambda  - \sum_s x_{t}^s \tfrac{\p\Lambda}{\p x^s}$.

Since the addition of the total derivative $-\tfrac{d}{dt}\Lambda$ to a Lagrangian does not change
the complete variation, the Euler-Lagrange equation with Lagrangian  \eqref{eq:EL1}
coincides with that substituted by $\tilde L= F(\bx,\bx_t)-\sum_s x_{t}^s \tfrac{\p \Lambda}{\p x^s}$.
We see that $ \tilde L$ is independent of $\bx_{tt}$ implying  the first claim of Proposition~\ref{prop:1}.
Next, since by our assumptions for a solution $\bx(t)$ any of its reparameterization  $\bx(\tau(t))$ is also a solution,
the Lagrangian $\tilde L$ is necessary homogeneous in $t$ of degree 1.

\section{PDE setup for the inverse problem} \label{S3}

Here we work with inhomogeneous ODE \eqref{ODE} and the corresponding Lagrangian $L$, which now
can be assumed of the first order. The variational problem
 \begin{equation}\label{Act}
\int L(x,\by,\dot{\by})\,dx\to\min
 \end{equation}
leads to the Euler-Lagrange equations
 \begin{equation}\label{EL}
\frac{\p L}{\p y^j}-\frac{d}{dx}\frac{\p L}{\p\dot{y}^j}=0,\quad 1\leq j\leq n,
 \end{equation}
where $\frac{d}{dx}=\p_x+\dot{y}^j\p_{y^j}+f^j\p_{\dot{y}^j}$ is the operator of total derivative.

This is an overdetermined 2nd order PDE system on a scalar function $L=L(x,\by,\dot{\by})$ and
it is equivalent to \eqref{ODE} if and only if:
Euler-Lagrange system \eqref{EL} vanishes modulo ODE system \eqref{ODE}
and the Jacobian matrix of $L$ is nondegenerate (to be able to express the ODE from the EL)
 $$
\det\left[\frac{\p^2 L}{\p{\dot{y}}^i\p{\dot{y}}^j}\right]_{i,j=1}^n\neq0.
 $$

Note that system \eqref{EL} is not of finite type, i.e.\ its solution is non-unique
(modulo divergences and  rescalings) and may be even not finitely parametric but contain arbitrary functions. 
Indeed, when $f^i=0$ the problem \eqref{Act} with straight lines as extremals has infinite-dimensional space of
solutions. These are the so-called Minkowski  Finsler metrics, given by translationally invariant Lagrangians
$L=L(\dot{\by})$. Clearly there is a functional freedom in choosing such a Finsler metric.

 % NB: This section is about $L$ not $\hat{L}$!
 % Indeed,  consider the so-called Minkowski  Finsler metrics, i.e., the Lagrangian  $F(x,x_t)$ such that it is homogeneous of degree $1$ in $x_t$ and such that it
 % does not depend on $x$. The extremals of $F$  generate flat projective structure and clearly there is a functional freedom in choosing such a Finsler metric.

\subsection{The fundamental system} \label{S31}

In \cite{Dav} Davis and in \cite{D} Douglas derived the following fundamental overdetermined system
on the symmetric nondegenerate matrix $\phi_{ij}=\frac{\p^2 L}{\p\dot{y}^i\p\dot{y}^j}=\phi_{ji}$:
 \begin{gather}
\frac{\p\phi_{ik}}{\p \dot{y}^j}=\frac{\p\phi_{jk}}{\p \dot{y}^i},\label{tg0}\\
\frac{d}{dx}\phi_{ij}+\frac12\frac{\p f^k}{\p\dot{y}^i}\phi_{kj}+\frac12\frac{\p f^k}{\p\dot{y}^j}\phi_{ki}=0,\label{tg1}\\
A_i^k\phi_{kj}=A_j^k\phi_{ki},\label{tg2}
 \end{gather}
where
 $$
A_i^j=\frac{d}{dx}\frac{\p f^j}{\p\dot{y}^i}-2\frac{\p f^j}{\p y^i}
-\frac12\frac{\p f^k}{\p\dot{y}^i}\frac{\p f^j}{\p\dot{y}^k}.
 $$
Note that $A$ is a (1,1)-tensor (or a field of operators, which is more obvious than in \cite{D},
when written in proper indices), so condition \eqref{tg2} means that $A$ is symmetric with respect to
metric $\phi$: $\phi(A\xi,\eta)=\phi(\xi,A\eta)$.
However $A$ is a given field and the unknown in this equation is $\phi$. Yet, there are many solutions
(depending on arbitrary functions).

To restrict those solutions further note that $A$, in general, is not integrable (its Nijenhuis tensor does not vanish),
and there are more constraints coming from \eqref{tg1}, and also \eqref{tg0}. Namely, passing from $A$ to
$A'=\frac{d}{dx}A-\frac12AJ-\frac12J^*A$, where $J^k_j=\frac{\p f^k}{\p\dot{y}^j}$ and $*$
is conjugation with respect to $\phi$ we get the equation $\phi(A'\xi,\eta)=\phi(\xi,A'\eta)$.
One can further iterate this recursive generation of constraints, and this is what is done in \cite{D}
for $n=2$. However, as we will see, for $n>2$ already the first iteration is generically sufficient.

\subsection{On reparametrizations} \label{S32}

If $\hat{L}=\hat{L}(\bx,\bx_t)$ is 1-homogeneous in velocity $\bx_t$, then the functional on curves in $M$
 $$
\bx(t)\mapsto \int \hat{L}(\bx(t),\bx_t(t))\,dt
 $$
is reparametrization invariant. In particular for a path $x^i=x^i(t)$ choosing $x^0=x$ instead of parameter
$t$ we obtain the integral in \eqref{Act}: indeed when $x^0=t$ we get $x^0_t=1$ and
 $$
L(x,y^j,\dot{y}^j)=\hat{L}(x,y^j,1,\dot{y}^j).
 $$

Conversely, given $L(x,y^j,\dot{y}^j)$ can be extended to a function 1-homogeneous in velocities
on $TM$ that is a cone over $STM$ as follows (for nonsymmetric $L$, i.e.\ if $L(-v)\neq L(v)$, $v\in ST_xM$,
one has to distinguish between $x^0_t<0$ and $x^0_t>0$ that may be not possible locally over domains in $M$,
but only microlocally on small domains $U\subset STM$):
 $$
\hat{L}(x^0,x^1,\dots,x^n,x^0_t,x^1_t,\dots,x^n_t)=
L\left(x^0,x^1,\dots,x^n,\frac{x^1_t}{x^0_t},\dots,\frac{x^n_t}{x^0_t}\right)\cdot x^0_t.
 $$

Recall that the condition for $\hat{L}$ to define a Finsler metric is the subadditivity
in velocities, which is equivalent (provided $\hat{L}$ smooth on $TM\setminus 0_M$)
to the strong convexity condition: for any $x\in M$ and $0\neq v\in T_xM$
the Hessian of $\hat{L}^2|_{T_xM}$ is positive definite at $v$.

Note that $\hat{L}^2$ is nondegenerate, i.e.\ $\det\text{Hess}\bigl(\hat{L}^2|_{CU}\bigr)\neq0$,
for the cone $CU\subset T_xM$ over an open dense subset $U\subset ST_xM$,
if $L$ is nonvanishing and nondegenerate:
 $$
\det\text{Hess}\bigl(\hat{L}^2|_{CU}\bigr)=2^{n+1}L^{n+2}\det\text{Hess}\bigl(L|_{U}\bigr)
 $$
(in general there are no relations between
nondegeneracy of $L$ and $L^2$.  Note also that, due to 1-homogeneuity,
$\det\text{Hess}\bigl(\hat{L}|_{T_xM}\bigr)\equiv0$).
We will call such $\hat{L}$ a pseudo-Finsler  metric (an example is $\sqrt{|g(v,v)|}$ for a Lorentzian metric $g$ on $M$). 
In this case equation $\{\hat{L}=1\}$ in $T_xM$
does not necessary define a convex but a nondegenerate (almost everywhere) hypersurface.

\section{Proof of Theorem \ref{Thm} } \label{S4}

In the case $f^i=f^i(x,\by)$ we have $A_i^j=-2\frac{\p f^j}{\p y^i}$ and the fundamental system and its
prolongation contian the following algebraic subsystem
 \begin{equation}\label{fy}
A_i^k\phi_{kj}=A_j^k\phi_{ki},\quad \bigl(\tfrac{d}{dx}A_i^k\bigr)\phi_{kj}=\bigl(\tfrac{d}{dx}A_j^k\bigr)\phi_{ki}.
 \end{equation}
This linear homogeneous system consists of $n(n-1)$ equations on $\frac12n(n+1)$ unknowns, and so
is determined for $n=3$ and overdetermined for $n>3$. We claim that generically it attains the maximal
rank, and hence the only solution is $\phi_{ij}=0$.

Note that for $n=2$ the system is underdetermined, hence as in \cite{D} one should add one more
linear equation\footnote{In the general case $\ff=\ff(x,\by,\dot{\by})$ the operator $A\mapsto A'$ of \S\ref{S31}
should be used instead of $\frac{d}{dx}$.}, namely $\bigl(\tfrac{d^2}{dx^2}A_i^k\bigr)\phi_{kj}=\bigl(\tfrac{d^2}{dx^2}A_j^k\bigr)\phi_{ki}$.
Then we get the $3\times3$ matrix of the system, which is generically nondegenerate, whence the same conclusion.
 % We consider henceforth the case $n>2$.

\subsection{Nonexistence of solutions to the inverse problem} \label{S41}

To prove the above claim for $n>2$ we first exhibit a system for which the maximal rank is attained.
This is given by
 \begin{equation}\label{fpar}
f^1=\sum_{k=1}^n(y^k)^2,\ f^2=(y^1)^2,\ f^3=(y^2)^2,\ \dots,\ f^{n-1}=(y^{n-2})^2,\ f^n=y^{n-1}.
 \end{equation}
The $n(n-1)\times\binom{n+1}{2}$ matrix $\mathbb{A}$ of system \eqref{fy}
(to obtain it write $\phi_{ki}$ into a column)
depends on $(x,\by,\dot{\by})$ and has maximal rank
$\binom{n+1}{2}$, for instance, at the point $x=0$, $y^j=\delta^j_n$, $\dot{y}^j=1$.
  (We omit this tedious verification.)
Since the rank is generally maximal and the
data are algebraic in 4-jets, the rank is generically maximal. Moreover, when we perturb the condition
$\frac{\p f^i}{\p\dot{y}^j}=0$ the matrix $\mathbb{A}$ changes, but still possesses the maximal rank,
and it will be maximal for generic 4-jets of the vector function $\ff$.
This implies the claim for $n>2$.

Let us give an alternative geometric, less computational agrument.
The matrices involved in \eqref{fy} have the form (note that $\mathbb{A}$ is not
this matrix below, but is easily derived from it):
 $$
\left(A\,\Bigl|\bigr.\,\tfrac{d}{dx}A\right)^t=-4\cdot\begin{pmatrix}
y^1 & y^1 & 0 & 0 & \dots & 0 & 0  &|& \dot{y}^1 & \dot{y}^1 & 0 & 0 & \dots & 0 & 0\\
y^2 & 0 & y^2 & 0 & \dots & 0 & 0 &|& \dot{y}^2 & 0 & \dot{y}^2 & 0 & \dots & 0 & 0\\
y^3 & 0 & 0 & y^3 & \dots & 0 & 0 &|& \dot{y}^3 & 0 & 0 & \dot{y}^3 & \dots & 0 & 0\\
\vdots & \vdots & \vdots & \vdots & \ddots & \vdots & \vdots &|&
\vdots & \vdots & \vdots & \vdots & \ddots & \vdots & \vdots\\
y^{n-2} & 0 & 0 & 0 & \dots & y^{n-2} & 0 &|& \dot{y}^{n-2} & 0 & 0 & 0 & \dots & \dot{y}^{n-2} & 0\\
y^{n-1} & 0 & 0 & 0 & \dots & 0 & \tfrac12 &|& \dot{y}^{n-1} & 0 & 0 & 0 & \dots & 0 & 0\\
y^n & 0 & 0 & 0 & \dots & 0 & 0 &|& \dot{y}^n & 0 & 0 & 0 & \dots & 0 & 0
\end{pmatrix}
 $$
For generic entries the blocks have different eigenvalues each and are mutually independent.
The quadraric form $\phi$ has eigenvectors of each block as orthogonal basis, but in
general two bases cannot be simultenously orthogonal for one metric (any signature).
This finished the proof.

For $n=2$ the same argument works with the same ODE system \eqref{fpar}.
In fact, this system for $n=2$ was already indicated by Douglas, and in \cite[formula (3.1)]{D}
the $3\times3$ matrix $\Delta$ is nondegenerate, implying $\phi_{ij}=0$ as the only solution.
Our observation extends his result without going into detailed analysis of solvability of the
fundamental system.

\subsection{Other approaches}\label{S42}

Let us consider one more example of nonexistence, namely
a higher-dimensional version of another system from \cite{D}:
 \begin{equation}\label{fpa2}
f^1=\sum_{k=1}^n(y^k)^2,\ f^2=0,\  \dots,\ f^{n-1}=0,\ f^n=0.
 \end{equation}
Then for the matrix $\mathbb{A}$ of system \eqref{fy} its $n\times n$ minor consisting of rows with numbers
$\bigl(1,\dots,n-1,\binom{n+1}2-1\bigr)$ and columns $(1,\dots,n)$ is equal to $(-2)\times$ the matrix
 $$
 \begin{pmatrix}
y^2 & y^1 & 0 & 0 & \cdots & 0 & 0\\
y^3 & 0 & y^1 & 0 & \cdots & 0 & 0\\
\vdots & \vdots & \vdots & \vdots & \ddots & \vdots & \vdots\\
y^{n-1} & 0 & 0 & 0 & \cdots & y^1 & 0\\
y^n & 0 & 0 & 0 & \cdots & 0 & y^1\\
\dot{y}^n & 0 & 0 & 0 & \cdots & 0 & \dot{y}^1
 \end{pmatrix}
\qquad\text{with }\quad \det=(y^1)^{n-2}(y^1\dot{y}^n-y^n\dot{y}^1)\not\equiv0,
 $$
while the columns $\bigl(n+1,\dots,\binom{n+1}2\bigr)$ vanish identically.
This implies that $\phi_{1i}=0$ and hence $\det(\phi_{ij})=0$.  Therefore ODE system \eqref{fpa2}
is not variational.

However this argument does not survive perturbation, as it belongs to a lower strata in branching
the compatibility analysis of system \eqref{tg0}-\eqref{tg2}. Complete analysis depending
on ranks of the arising matrices was performed for $n=2$ in \cite{D}. However the number of branches
grows rapidly with $n>2$ and it would be unreasonable to expect a complete answer due to complexity.

One can exploit the idea of \cite{KM2} to find the number of independent solutions (dimension) when
system \eqref{tg0}-\eqref{tg2} is of finite type. Namely, its prolongation, obtained by differentiation of
all equations of the system to a sufficiently large order $N$ at a particular point $z\in U$,
stabilizes the solution space, given by $(N+1)$-jet of $\phi_{ij}$ at $z$.
In practice this procedure allows to effectively decide solvability of the system.

\section{Submaximally symmetric structures are variational}\label{S5}

In this section we discuss several examples, where we can resolve the fundamental system for the inverse problem.
Namely we consider path structures admitting infinitesimal symmetries, i.e., local diffeomorphisms preserving the structure.
A flat structure on a manifold $M$ of dimension $n+1$ has maximally symmetry dimension $n^2+4n+3$ and is is variational.

The next submaximal symmetry dimension is $n^2+5$ for $n>1$; let us specify submaximal symmetry 
depending on the type of (nonzero) harmonic curvature,
namely Fels torsion $T$ or Fels curvature $S$, see \cite[\S5.3-5.4]{KL}. 
In the zero curvature module ($S=0$) we get projective geometry, 
and in the torsion-free module ($T=0$) we het general path geometry (Segr\'e branch; non-projective).
We consider those in turn.

\subsection{The Egorov projective structure} \label{S51}

This structure is originally \cite{E} given by the nonzero Christoffel coefficients
$\Gamma_{23}^1=\Gamma_{32}^1=x^2$ on $M=\mathbb{R}^{n+1}(\bx)$.
The corresponding inhomogeneous system \eqref{ODE} has the following form:
 \begin{equation}\label{Eg}
\ddot{y}^j=2y^1\dot{y}^1\dot{y}^2\dot{y}^j,\qquad 1\leq j\leq n.
 \end{equation}
This structure has maximal symmetry dimension $n^2+4$ among all nonflat projective structures \cite{E,T} and up to local diffeomorphism it is 
unique such \cite{T}; it is non-metrizable by \cite{KM1}, i.e.\ there is no Levi-Civita connection in its projective class.

Surprisingly, the structure  is variational, at least micro-locally:

 \begin{proposition}\label{Prop}
There exists a Lagrangian function $L$ defined for almost all velocities, which generates the Egorov projective structure.
 \end{proposition}

To see this note first following \cite[remark after Theorem 3]{BLP}
that equation \eqref{Eg} is linearizable, namely it transforms under a point transformation to the ODE
 \begin{equation}\label{EgL}
\ddot{y}^1=y^2,\ \ddot{y}^2=0,\ \dots,\ \ddot{y}^n=0.
 \end{equation}
We will treat therefore this system. It is precisely of the kind considered at the beginning of this section.
Thus considering system \eqref{tg0}-\eqref{tg2} for this choice of $f^i$ we find a Lagrangian
 $$
L=(\dot{y}^1-x\,y^2)\dot{y}^2+\sum_3^n(\dot{y}^i)^2
 $$
with extremals given by \eqref{EgL}. The corresponding 1-homogeneous Lagrangian is
 $$
\hat{L}=\left(\frac{x^1_t}{x^0_t}-x^0x^2\right)x^2_t+\frac{(x^3_t)^2+\dots+(x^n_t)}{x^0_t}.
 $$
Its extremal curves satisfy the (underdetermined) ODE with the same paths as \eqref{EgL}:
 $$
x^0_tx^1_{tt}-x^1_tx^0_{tt}=(x^0_t)^3x^2,\ x^0_tx^j_{tt}-x^j_tx^0_{tt}=0\quad (1<j\leq n).
 $$
Thus the Egorov structure is  variational.

 \begin{remark} 
Lagrangians of the form $\hat{L}=\frac{g(\bx_t,\bx_t)}{\alpha(\bx_t)}$ for a Riemannian or pseudo-Riemannian metric $g$ and for a
1-form $\alpha$ are called Kropina (pseudo-Finsler) metrics. Kropina metrics were also considered in the framework of mathematical relatively, 
see e.g.\ E. Caponio et al \cite{Caponio1}.
Kropina metrics are not defined on the vectors $\bx_t$ lying in the kernel of $\alpha$. Note that in our case the form $\alpha$ is closed 
so its extremals define a projective structure by \cite[Corollary 3.6]{Cheng}.
 \end{remark}

In this context, the following question is natural: does there exist  a strictly convex Finsler metric (without singularities and defined on the whole slit tangent bundle)  whose geodesics are curves of the Egorov projective structure? The next proposition answers this question negatively:

 \begin{proposition}
The path structure   given by \eqref{EgL}, and hence by \eqref{Eg}, is not Finsler metrizable. % does not come from any Finsler metric.
 \end{proposition}

Indeed, in this case we can obtain the general solution of the fundamental system, which due to
a very simple form $A_i^j=-2\delta_1^j\delta^2_i$, is as follows:
 $$
\phi_{11}=\phi_{13}=\dots=\phi_{1n}=0,\  \phi_{12}\neq0,\ \tfrac{d}{dx}\phi_{ij}=0,
 $$
where $\frac{d}{dx}=\p_x+\dot{y}^i\p_{y^i}+y^2\p_{\dot{y}^1}$.
This implies the form (we omit dependence of $\psi_0,\psi_1$ on $\bx$, indicating only velocity $\bx_t$)
of the homogeneous Lagrangian:
 $$
\hat{L}= \psi_0\left(\frac{x^2_t}{x^0_t},\dots,\frac{x^n_t}{x^0_t}\right)x^0_t+
\psi_1\left(\frac{x^2_t}{x^0_t},\dots,\frac{x^n_t}{x^0_t}\right)x^1_t.
 $$
One can easily see that for any choice of $\psi_0,\psi_1$ the function $\hat{L}^2$ is not convex.

\subsection{Submaximally symmetric path structure}\label{S52}

The maximally symmetric nonflat path structure has dimension of the symmetry algebra $n^2+5$, see \cite{KT}.
Uniqueness of such a structure has been recently established in \cite{T}.
For $n+1=3$ this structure was discussed in \cite{CDT} in relation to self-dual gravity,
the corresponding spacetime\footnote{This Plebanski type metric has coordinate expression
$g=dx\,dw+dy\,dz-y^2dw^2$.} is Ricci flat of Petrov type N.
The ODE system generating this metric via the twistor correspondence is
 \begin{equation}\label{submax}
\ddot{y}^1=(\dot{y}^2)^3,\ \ddot{y}^2=0,\ \dots,\ \ddot{y}^n=0.
 \end{equation}
The fundamental system for the inverse problem is solvable; one solution is given by
 $$
L=\left(\sqrt{\pi}\,\dot{y}^1\op{erf}\left(\frac{\dot{y}^1}{(\dot{y}^2)^{3/2}}\right)
+(\dot{y}^2)^{3/2} \exp\left(-\frac{(\dot{y}^1)^2}{(\dot{y}^2)^3}\right) + \dot{y}^1\right) e^{2y^1}
+ \sum_3^n (\dot{y}^i)^2.
 $$
Thus the path structure  \eqref{submax} is variational.
The corresponding 1-homogeneous Lagrangian $\hat{L}$ can be derived straightforwardly.

\subsection{More examples}\label{S53}

Another notable path structure  is given by a family of distinguished curves of the trivial scalar
ODE, encoded as the flat $A_2/P_{1,2}$ homogeneous geometry \cite{CZ}. The distinguished curves
transversal to the contact structure on $J^1(\mathbb{R}^1)$ are given by a pair of differential equations
on unknowns $y^1(x),y^2(x)$ (cf. an equivalent form in \cite[\S7.2]{CDT}):
 \begin{equation}\label{distinguished}
\ddot{y}^1=\frac{2(\dot{y}^1)^2}{y^1-\dot{y}^2},\ \ddot{y}^2=0.
 \end{equation}
This ODEs system is also related to anti self-dual conformal metrics, namely it 
generates via the twistor correspondence an Einstein metric of constant negative scalar curvature \cite{CDT}.

The fundamental system passes the compatibility test (as discussed at the end of \S\ref{S42}),
so from the Cartan-K\"ahler theorem it follows that it possesses solutions with any admissible Cauchy initial data;
in particular, we conclude that system \eqref{distinguished} is variational. Indeed, for
 $$
L=\frac{\dot{y}^1}{\dot{y}^2-y^1}
 $$
extremals are exactly the paths given by \eqref{distinguished}. The corresponding 1-homogeneous
Lagrangian $\hat{L}$ can be derived straightforwardly.

 % Our first attempts to solve this overdetermined system via Maple failed, as \textsf{pdsolve}
 % returned empty output. Then we used a trick for $n=2$ by presenting unknown Lagrangian as
 % a rational function, where the numerator was quadratic in velocities, while denominator was
 % linear, whose vanishing corresponds to the contact structure (this idea was borrowed from
 % \cite{MMM} in the context of CR geometry). This yields the following Lagrangian
 
 \begin{remark}
An elliptic version of this example consist of chains in (not necessary flat) CR geometry. 
It was proven in \cite{Cheng} that in any dimension the path geometry of chains is variational, 
with the Lagrangian being a Kropina metric.
 \end{remark}

\section{Conclusion}

The inverse variational problem for nonautonomous ODE systems \eqref{ODE} has attracted a lot of interest 
in the literature; several criteria for variationability were obtained.
We have shown that a generic path structure in dimension  $n+1\ge 3$ is not variational.
The proof is done in terms of jets. % and works also micro-locally.
Our methods allow to derive a proper subanalytic subset $\Sigma\subset J^k$ such that (regular) variational structures
given as \eqref{ODE} are subject to the constraints $j^\ell\ff(U)\subset\Sigma$, $U\subset STM$.
 % As demonstrated in \S\ref{S51}, this allows in certain cases to reconstruct the Lagrangian by a path structure.

In particular if a path structure comes from experimental observations and should be variational by
physical reasons, our methods may help to confirm correctness of the experiment;
and also find a variational structure that (in some sense) is closest to the experimental data. 
We leave aside a related question on optimal regularity $C^k$, where our results hold. 

A corollary of our main theorem implies that a generic projective structure is not metrizable in the class of
Riemannian or pseudo-Riemannian metrics.
This result was expected: indeed, the metrization problem can be reduced to an overdetermined system of PDEs
of finite type (see e.g. M.  Eastwood et al \cite{Eastwood}).
Nevertheless, this result was formally established only in dimension 2 (R. Bryant et al \cite{Bryant})
and in dimension 3 (M. Dunajski et  al \cite{Dunajski}); we proved  it in any dimension.

We also demonstrated that the Egorov projective structure is variational in any dimension $n+1\ge3$
by exhibiting a Kropina type pseudo-Finsler metric. For $n=2$ this could be obtained from the
results of Douglas \cite{D}; in \cite{Anderson} another Lagrangian was derived for 
ODE \eqref{EgL} though without any relation to the Egorov structure.
By our previous work \cite{KM1} it is not metrizable
in the pseudi-Riemannian setting. In this work we proved it is not metrizable in the 
Finsler setting. We also demonstrated variationability of some other notable path geometries with
many infinitesimal symmetries.


\begin{thebibliography}{100}

\bibitem{Berck}
J.-C. Alvarez-Paiva, G. Berck, {\it Finsler surfaces with prescribed geodesics}, arXiv:1002.0243 (2010).

\bibitem{Anderson}
I. Anderson, G.  Thompson, {\it The inverse problem of the calculus of variations for ordinary differential equations},
Mem. Amer. Math. Soc. {\bf 98}, no. 473 (1992).

\bibitem{Beltrami}
E. Beltrami, {\it Resoluzione del problema: riportari i punti di una superficie sopra un piano in modo che le linee
geodetische vengano rappresentante da linee rette}, Annali di Matematica % pura et applicata
\textbf{1}, no. 7, 185--204 (1865).

\bibitem{BLP}
V.\ Boyko, O.\ Lokaziuk, R.\ Popovych, {\it Admissible transformations and Lie symmetries
of linear systems of second-order ordinary differential equations}, arXiv:2105.05139 (2021).

\bibitem{Bryant0}
R. Bryant, G. Manno, V. S. Matveev, {\it A solution of a problem of Sophus Lie: normal forms of two-dimensional
metrics admitting two projective vector fields}, Math. Ann. {\bf 340}, no. 2, 437--463 (2008).
		
\bibitem{Bryant}
R. Bryant, M. Dunajski, M.  Eastwood, {\it Metrisability of two-dimensional projective structures},
J. Differential Geom. { \bf  83}, no. 3, 465--499 (2009).

\bibitem{burns} 
K.\ Burns, V. Matveev, {\it Open problems and questions about geodesics}, 
Ergodic Theory Dynam. Systems {\bf 41}, no. 3, 641--684 (2021). 

\bibitem{Busemann0}
H. Busemann, {\it Two-dimensional metric spaces with prescribed geodesics},
Ann.\ Math.\ (2) {\bf 40}, no.1, 129--140 (1939).
		
\bibitem{Busemann_book}
H. Busemann, {\it The geometry of geodesics}, Academic Press Inc, New York (1955).
		
\bibitem{CZ}
A.\ \v{C}ap, V.\ \v{Z}\'adnik, {\it Contact projective structures and chains}, Geom.\ Dedicata {\bf 146}, 67--83 (2010).

\bibitem{Caponio}
E. Caponio, M.A. Javaloyes, M. Sanchez, {\it On the interplay between Lorentzian causality and Finsler metrics
of Randers type}, Rev. Mat. Iberoam. {\bf 27}, no. 3, 919--952 (2011).

\bibitem{Caponio1}
E. Caponio, M.A. Javaloyes, M. Sanchez, {\it Wind Finslerian structures: from Zermelo's navigation to the
causality of spacetimes}, to appear in Memoirs AMS; arXiv:1407.5494 (2014).

\bibitem{C}
E.\ Cartan, {\it Sur les vari\'et\'es \`a connexion projective}, Bulletin de la S.M.F.\ {\bf 52}, 205-241 (1924).

\bibitem{CDT}
S.\ Casey, M.\ Dunajski, P.\ Tod, {\it Twistor geometry of a pair of second order ODEs},
Comm.\ Math.\ Phys.\ {\bf 321}, 681--701 (2013).
  	
\bibitem{Cheng} J.-H. Cheng, T.  Marugame, V. S.  Matveev, R.  Montgomery,
{\it Chains in CR geometry as geodesics of a Kropina metric, }
Adv.\ Math.\ {\bf 350}, 973--999 (2019).
		
\bibitem{Dar}
G.\ Darboux,  {\it Le\c{c}ons sur la th\'eorie g\'en\'erale des surfaces},
Vol. III, \S 604--605, Gauthier-Villars, Paris (1894).

\bibitem{Dav}
D.\ R.\ Davis, {\it The inverse problem of the calculus of variations in a space of $n+1$ dimensions},
Bull. Amer.\ Math.\ Soc.\ {\bf 35}, 371--380 (1929).

\bibitem{Dini}
U.\ Dini, {\it Sopra un problema che si presenta nella teoria generale delle rappresentazioni geografice di una superficie su un'altra},
Ann.\ di Math.\ Ser.2 {\bf 3}, 269--293 (1869).

\bibitem{Do1}
T.\ Do, G.\ Prince, {\it New progress in the inverse problem in the calculus of variations,}  Differential Geom. Appl. {\bf 45}, 148--179 (2016).

\bibitem{Do}
T.\ Do, G.\ Prince, {\it The inverse problem in the calculus of variations: new developments},
Commun. Math. {\bf  29}, no. 1, 131--149 (2021).

\bibitem{D}
J.\ Douglas, {\it Solution of the inverse problem of the calculus of variations},
Transactions Amer.\ Math.\ Soc.\ {\bf 50}, 71--128 (1941).

\bibitem{Dunajski}
M.\ Dunajski, M.\ Eastwood, {\it Metrisability of three-dimensional path geometries},
Eur.\ J.\ Math.\ {\bf 2}, 809--834 (2016).

\bibitem{Eastwood}
M.\ Eastwood, V.\ Matveev, {\it Metric connections in projective differential geometry}, in:
{\it Symmetries and overdetermined systems of partial differential equations,} 339--350,
IMA Vol. Math. Appl., {\bf 144}, Springer, New York (2008).

\bibitem{E}
I.\,P.\ Egorov, {\it Collineations of projectively connected spaces},
Doklady Akad. Nauk SSSR {\bf 80}, 709--712 (1951).

\bibitem{EPS}
J.\  Ehlers, F.\  Pirani, A.\  Schild, {\it The geometry of free fall and light propagation}, in
{\it General Relativity, Papers in Honour of J.L.\ Synge}, ed.\ L.\ O'Raifertaigh,
Oxford:  Clarendon Press pp.~63--84  (1972);
Republished in {\it General Relativity and Gravity}, {\bf 44}, 1587--1609 (2012).

\bibitem{GM}
J.\ Grifone, Z.\ Muzsnay, {\it Variational principles for second-order differential equations. Application of the Spencer theory to characterize
variational sprays}, World Scientific Pub (2000).

\bibitem{H} 
H.\ Helmholtz, {\it  \"Uber der physikalische Bedeutung des Princips der kleinsten Wirkung}, J.
Reine Angew. Math. {\bf 100}, 137--166 (1887).

\bibitem{Kr}
B. Kruglikov, {\it Point classification of second order ODEs: Tresse classification revisited and beyond\/}
(with an appendix by B.Kruglikov and V.Lychagin), in {\it  Abel Symp. {\bf 5},
Differential equations: geometry, symmetries and integrability,}  199--221, Springer, Berlin (2009). 

\bibitem{KL}
B. Kruglikov, V. Lychagin, {\it Geometry of Differential equations,}  in: {\it Handbook of Global Analysis, }
Eds: D.Krupka, D.Saunders, Elsevier, 725--772 (2008). 

\bibitem{KM1}
B.\ Kruglikov, V.\ Matveev, {\it Submaximal metric projective and metric affine structures },
Diff.\ Geom.\ Appl. {\bf 33}, 70--80 (2014).

\bibitem{KM2}
B.\ Kruglikov, V.\ Matveev, {\it Nonexistence of an integral of the 6th degree in momenta for the Zipoy-Voorhees metric},
Physical Review D {\bf 85}, 124057 (2012).

 % \bibitem{KM3}
 % B.\ Kruglikov, V.\ Matveev, {\it Strictly nonproportional geodesically equivalent metrics have $h_{top}(g)=0$},
 % Ergod.\ Th.\ & Dynam.\ Sys.\ {\bf 26}, 219-–245 (2006).

\bibitem{KT}
B.\ Kruglikov, D.\ The, {\it The gap phenomenon in parabolic geometries},
Journal f\"ur die Reine und Angew.\ Math.\
{\bf 723}, 153--215 (2017).

\bibitem{Krupka}
D. Krupka, {\it The Vainberg-Tonti Lagrangian and the Euler-Lagrange mapping},
Differential Geometric Methods in Mechanics and Field Theory, Eds: F. Cantrijn, B. Langerock,
% Volume in Honour of W. Sarlet,
Gent Academia Press, 81--90 (2007).

\bibitem{Lang}
J. Lang, {\it Finsler metrics on surfaces admitting three projective vector fields},
Differential Geom. Appl. {\bf 69}, 101590 (2020).

 \bibitem{LC}
T.\  Levi-Civita, {\it Sulle trasformazioni delle equazioni dinamiche}, Ann. di Mat., serie $2^a$,
{\bf 24}, 255--300 (1896).

\bibitem{Lie}
S.\  Lie, {\it Untersuchungen \"uber geod\"atische Kurven}, Math. Ann. \textbf{20} (1882);
Sophus Lie Gesammelte Abhandlungen, Band {\bf 2}, erster Teil, 267--374. Teubner, Leipzig (1935).

\bibitem{Read}
N.\ Linnemann, J.\ Read, {\it Constructive Axiomatics in Spacetime Physics Part I: Walkthrough to the
Ehlers-Pirani-Schild Axiomatisation}, arXiv:2112.14063 (2021).

\bibitem{matveev}
V.\ S.\ Matveev, {\it Geodesically equivalent metrics in general relativity},
J. Geom. Phys. {\bf 62}, no. 3, 675--691 (2012).
 
\bibitem{alone} 
V. S. Matveev, {\it Two-dimensional  metrics admitting  precisely one projective vector field,}  (with an appendix by A. Bolsinov V. S. Matveev  and G. Pucacco), Math. Ann. {\bf 352}, no. 4, 865--909 (2012).  

\bibitem{matveev2}
V.\,S.\ Matveev, {\it On projective equivalence and pointwise projective relation of Randers metrics},
Internat.\ J.\ Math.\ {\bf 23}, no. 9, 1250093, 14 pp (2012).

\bibitem{Scholz}
V.\,S.\ Matveev, E.\ Scholz, {\it Light cone and Weyl compatibility of conformal and projective structures},
Gen. Relativity Gravitation {\bf 52}, no. 7, Paper No. 66, 9 pp (2020).

 % \bibitem{MT}
 % V.\,S.\ Matveev, P.\,J.\ Topalov, {\it Trajectory equivalence and corresponding integrals},
 % Regular and Chaotic Dynamics {\bf 3} (2), 30--45 (1998).

\bibitem{Pfeiffer}
Ch.\ Pfeifer, {\it Finsler spacetime geometry in physics}, International Journal of Geometric Methods in Modern Physics
{\bf 16}, Supp. 2,  1941004 (2019).

\bibitem{Schur}
F. Schur, {\it Ueber die Deformation der R\"aume constanten Riemann'schen Kr\"ummungsmaasses},
Math. Ann. {\bf 27}, no. 2, 163--176 (1886).

\bibitem{Sonin}
N.\,Ya.\ Sonin, {\it On the definition of maximal and minimal properties} (in Russian), Warsaw Univ. Izvestiya
{\bf 1-2}, 1--68 (1886).

\bibitem{Tabachnikov} 
S. Tabachnikov, {\it Remarks on magnetic flows and magnetic billiards, Finsler metrics and a magnetic analog of Hilbert's fourth problem}, in: {\it Modern dynamical systems and applications}, 233--250, Cambridge Univ. Press, Cambridge (2004).   

\bibitem{T}
D.\ The, {\it On uniqueness of submaximally symmetric parabolic geometries}, arXiv:2107.10500 (2021).

\bibitem{Tresse}
A.\ Tresse, {\it D\'etermination des invariants ponctuels de l\'equation diff\'erentielle ordinaire du second ordre $y'' = \omega (x,y,y')$},
Leipzig.\ {\bf 87} S.\,gr.\ $8^\circ$ (1896).

\bibitem{Veblen}
O.\ Veblen, T.\,Y.\ Thomas, {\it The geometry of paths}, Trans.\ Amer.\ Math.\ Soc.\ {\bf 25}, no. 4, 551--608 (1923).

\bibitem{Weyl}
H.\ Weyl, {\it Zur Infinitisimalgeometrie: Einordnung der projektiven und der  konformen Auffasung},
Nachrichten von der K. Gesellschaft der Wissenschaften zu G\"ottingen, Mathematisch-Physikalische Klasse, 1921;
``Selecta Hermann Weyl'', Birkh\"auser Verlag, Basel und Stuttgart (1956).

\bibitem{Weyl1}
H.\ Weyl, {\it Mathematische {A}nalyse des {R}aumproblems. Vorlesungen gehalten in Barcelona und Madrid},
Berlin etc (1923); Springer. Nachdruck Darmstadt: Wissenschaftliche Buchgesellschaft (1963).

\end{thebibliography}
\end{document}